\input amstex
\documentstyle{amsppt}
\magnification=\magstep1                        
\hsize6.5truein\vsize8.9truein                  
\NoRunningHeads
\loadeusm

\magnification=\magstep1                        
\hsize6.5truein\vsize8.9truein                  
\NoRunningHeads
\loadeusm

\document
\topmatter

\title  
A note on Barker polynomials
\endtitle

\rightheadtext{A note on Barker polynomials}

\author Peter Borwein and Tam\'as Erd\'elyi
\endauthor

\address  Department of Mathematics and Statistics,
Simon Fraser University, Burnaby, B.C., Canada V5A 1S6 (P. Borwein) \endaddress

\address Department of Mathematics, Texas A\&M University,
College Station, Texas 77843 (T. Erd\'elyi) \endaddress

\email pborwein\@sfu.ca (P. Borwein) 
and terdelyi\@math.tamu.edu (T. Erd\'elyi)
\endemail

\endthanks


\abstract
We call the polynomial $P_{n-1}(x) = \sum_{j=1}^n{a_jz^{n-j}}$ 
a Barker polynomial of degree $n-1$ if each $a_j \in \{-1,1\}$ and
$$T_{n}(z) = P_{n-1}(z)P_{n-1}(1/z) = n + \sum_{k=1}^{n-1}{c_k(z^k + z^{-k})}\,, \qquad |c_k| \leq 1 \,.$$ 
Properties of Barker polynomials were studied by Turyn and Storer 
thoroughly in the early sixties, and by Saffari in the late eighties.  In the last few years 
P. Borwein  and his collaborators revived interest in the study of Barker polynomials (Barker 
codes, Barker sequences). In this paper we give a new proof of the fact that there is no 
Barker polynomial of even degree greater than $12$, and hence Barker sequences of odd length 
greater than $13$ do not exist. This is intimately tied to irreducibility questions and proved 
as a consequence of the following new result.

\proclaim{Theorem}
If $n : = 2m+1 > 13$ and  
$$Q_{4m}(z) = (2m+1)z^{2m} + \sum_{j=1}^{2m}{b_j(z^{2m-j} + z^{2m+j})}$$
where each $b_{j} \in \{-1,0,1\}$ for even values of $j$, each $b_j$ is
an integer divisible by $4$ for odd values of $j$, then
there is no polynomial $P_{2m} \in {\Cal L}_{2m}$ such that
$Q_{4m}(z) = P_{2m}(z)P_{2m}^*(z)$, where
$$P_{2m}^*(z):=z^{2m}P_{2m}(1/z)\,,$$
and ${\Cal L}_{2m}$ denotes the collection of all polynomials of degree $2m$ 
with each of their coefficients in $\{-1,1\}$.
\endproclaim

\noindent A clever usage of Newton's identities plays a central role in our elegant proof. 

\endabstract

\endtopmatter

\document

Following the pioneering work of Turyn and Storer [TS-61], Turyn [T-63], [T-65], and 
Saffari [S-90], in the last few years P. Borwein and his collaborators revived interest 
in the study of Barker polynomials (Barker codes, Barker sequences). See 
[B-02], [BCJ-12], [BM-08], [H-09], [M-09], and [E-11], for example. 

Following the paper of Turyn and Storer [TS-61] we call the polynomial 
$$P_{n-1}(x) = \sum_{j=1}^n{a_jz^{n-j}} = 
a_1z^{n-1} + a_2z^{n-2} + \cdots + a_{n-1}z + a_n \,, \qquad a_j \in \{-1,1\}\,,$$
a Barker polynomial of degree $n-1$ if
$$T_{n}(z) = P_{n-1}(z)P_{n-1}(1/z) = n + \sum_{k=1}^{n-1}{c_k(z^k + z^{-k})}\,, 
\qquad |c_k| \leq 1 \,.$$
Let $c_0 := n$. Since 
$$c_k = \sum_{j=1}^{n-k}{a_ja_{j+k}}\,, \qquad k=0,1, \ldots n-1\,,$$
and each term ${a_ja_{j+k}}$ is in $\{-1,1\}$, it follows that
$$\prod_{j=1}^{n-k}{(a_ja_{j+k})} = (-1)^{(n-k-c_k)/2}\,, 
\qquad k=0,1,\ldots, n-1\,. \tag 1$$
Multiplying two consecutive equations of the above form gives
$$a_{k+1}a_{n-k} = (-1)^{n-k-(1+c_k+c_{k+1})/2}\,, \qquad k=0,1,\ldots,n-1\,.$$ 
Observe that
$$c_k+c_{n-k} = \sum_{j=1}^{n}{a_ja_{j+k}}\,, \qquad k=1,2,\ldots,n-1\,,$$
where the second index $j+k$ is taken $\text {\rm modulo\,\,} n$.
Then (1) yields
$$\prod_{j=1}^{n}{(a_ja_{j+k})} = 1 = (-1)^{(n-c_k-c_{n-k})/2}\,,
\qquad k=0,1,\ldots, n-1\,,$$
hence $c_k+c_{n-k}=n \enskip (\text {\rm mod\,\,} 4)$.

Now assume that $n=2m+1$ is odd,  
$$|c_k| \leq 1\,, \qquad k=2,4,\ldots, 2m=n-1\,,$$ 
and each $c_k$ is an integer divisible by $4$ whenever $k=1,3,\ldots,2m-1=n-2$.  
Observe that $c_k+c_{n-k}=n \enskip (\text {\rm mod\,\,} 4)$ implies 
$$c_{2j-1}=0 \enskip (\text {\rm mod\,\,} 4) \qquad \text {\rm and} 
\qquad c_{2j}=(-1)^{m}\,, \qquad j=1,2,\ldots,m\,.$$
The formula for $a_{k+1}a_{n-k}$ becomes
$$a_{k+1}a_{n-k}=(-1)^{m+k}\,, \qquad k=0,1,\ldots,n-1\,. \tag 2$$  
Hence
$$c_{2j} = \sum_{i=1}^{n-2j}{a_ia_{i+2j}}$$
can be rewritten as
$$(-1)^{m} = \sum_{i=1}^{n-2j}{a_ia_{n+1-i-2j}(-1)^{m+i+1}}\,, \qquad j=1,2,\ldots\,m.$$
Now let $n-2j=2k+1, k=0,1,\ldots,m-1$. Then
$$\split 1 = & \sum_{i=1}^{2k+1}{a_ia_{2k+2-i}(-1)^{i+1}} \cr 
= & 2\,\sum_{i=1}^k{a_ia_{2k+2-i}(-1)^{i+1}} + a_{k+1}^2(-1)^{k+2}\,, 
\qquad k=0,1,\ldots,m-1\,. \cr \endsplit $$
Hence
$$\frac{1 + (-1)^{k+1}}{2} = \sum_{i=1}^k{a_ia_{2k+2-i}(-1)^{i+1}}\,,
\qquad k=0,1,\ldots, m-1\,. \tag 3$$
From this we can easily deduce
$$\prod_{i=1}^k{(a_ia_{2k+2-i})}=1\,,$$
and hence
$$\prod_{i=1}^{2k+1}{a_i}=a_{k+1}\,, \qquad k=0,1,\ldots,m-1\,.$$
Multiplying two consecutive equations of the above form gives
$$a_ka_{k+1} = a_{2k}a_{2k+1}, \qquad k=1,2,\ldots, m-1\,. \tag 4$$

Let ${\Cal L}_n$ be the collection of all polynomials of degree $n$ with each of
their coefficients in $\{-1,1\}$.

\proclaim{Theorem 1}
If $n : = 2m+1 > 13$ and  
$$Q_{4m}(z) = (2m+1)z^{2m} + \sum_{j=1}^{2m}{b_j(z^{2m-j} + z^{2m+j})}$$
where each $b_{j} \in \{-1,0,1\}$ for even values of $j$, each $b_j$ is
an integer divisible by $4$ for odd values of $j$, then
there is no polynomial $P_{2m} \in {\Cal L}_{2m}$ such that
$Q_{4m}(z) = P_{2m}(z)P_{2m}^*(z)$, where
$$P_{2m}^*(z):=z^{2m}P_{2m}(1/z)\,.$$
\endproclaim

Observe that if $Q_{4m}(z)$ satisfies the conditions of the theorem and there 
is a polynomial $P_{2m} \in {\Cal L}_{2m}$ such that
$Q_{4m}(z) = P_{2m}(z)P_{2m}^*(z)$, where 
$$P_{2m}^*(z):=z^{2m}P_{2m}(1/z)\,,$$
then the introductory remarks imply that 
$$Q_{4m}(z) = P_{2m}(z)P_{2m}^*(z) = (-1)^mP_{2m}(z)P_{2m}(-z)\,,$$
that is the polynomial $Q_{4m}$ is even, hence, in fact, 
we have that $b_j=0$ for each odd values of $j$. Moreover 
$$b_j = (-1)^m\,, \qquad j=2,4,\ldots,2m\,.$$
Therefore to prove Theorem 1 it is sufficient to prove only the following much 
less general looking result. 

\proclaim{Theorem 2}
Suppose $n : = 2m+1 > 13$ and
$$Q_{4m}(z) = (2m+1)z^{2m} + (-1)^m\sum_{j=1}^m{(z^{2m-2j} + z^{2m+2j})}\,.$$
There is no polynomial $P_{2m} \in {\Cal L}_{2m}$ such that 
$Q_{4m}(z) = P_{2m}(z)P_{2m}^*(z)$, where 
$$P_{2m}^*(z):=z^{2m}P_{2m}(1/z)\,.$$
As a consequence, Barker polynomials $P_{2m}$ of degree greater than $12$ 
do not exist.
\endproclaim

Before presenting the proof of Theorem 2 we need to introduce some 
notation and to prove two lemmas. 

Suppose $P_{2m} \in {\Cal L}_{2m}$ is of the form
$$P_{2m}(z) = \sum_{j=1}^{2m+1}{a_jz^{2m+1-j}}\,, \qquad a_j \in \{-1,1\}\,, 
\quad j=1,2,\ldots,2m+1\,,$$
and assume that
$$Q_{4m}(z) = P_{2m}(z)P_{2m}^*(z)\,.$$
Note that (2) implies that $P_{2m}^*(z)=(-1)^mP_{2m}(-z)$.
Without loss of generality we may assume that $a_1=1$.
We introduce the coefficients $c_j$ of $Q_{4m}$ by
$$Q_{4m}(z) = \sum_{j=1}^{4m+1}{c_{4m+2-j}z^{j-1}}\,,$$
that is, $c_{2m+1} = 2m+1$ and
$$c_1 = c_3 = \cdots = c_{2m-1} = c_{2m+3} = c_{2m+5} = \cdots = c_{4m+1} = (-1)^m\,.$$
Let 
$$P_{2m}(z) = \prod_{j=1}^{2m}{(z-\alpha_j)}\,, \qquad \alpha_j \in {\Bbb C}\,.$$   
Then
$$P_{2m}(-z) = \prod_{j=1}^{2m}{(z + \alpha_j)}\,, \qquad \alpha_j \in {\Bbb C}\,,$$
and
$$Q_{4m}(z) = (-1)^m \prod_{j=1}^{2m}{(z - \alpha_j)(z + \alpha_j)}\,, 
\qquad \alpha_j \in {\Bbb C}\,.$$
Associated with a nonnegative integer $\mu$ let 
$$S_{\mu} = S_{\mu}(Q_{4m}) := \sum_{j=1}^{2m}{\alpha_j^{\mu} + (-\alpha_j)^{\mu}}\,,$$
and
$$s_{\mu} = s_{\mu}(P_{2m}) := \sum_{j=1}^{2m}{\alpha_j^{\mu}}\,.$$
Observe that $S_{\mu} = 0$ when $\mu$ is odd.
Newton's identities (see page 5 in [BE-95], for instance) give
$$c_1S_\mu + c_2S_{\mu-1} + \cdots + c_{\mu}S_1 + \mu c_{\mu+1} = 0\,, 
\qquad \mu =1,2,\ldots, 4m\,,$$
that is,
$$S_{\mu} + S_{\mu-2} + S_{\mu-4} + \cdots + S_2 = -\mu\,, 
\qquad \mu = 2,4,\ldots, 2m-2\,.$$
We conclude that
$$S_2 = S_4 = S_6 = \cdots = S_{2m-2} = -2\,,$$
and since $2s_{2k} = S_{2k}$, we have 
$$s_{2k} = -1\,, \qquad k=1,2,\ldots, m-1\,.$$
Let
$$a_1 = a_2 = \cdots = a_p =1, \quad a_{p+1} = -1\,.$$
Then (4) implies that $p$ is odd and, considering $P_{2m}^*$ rather than 
$P_{2m}$ if necessary, without loss of generality we may assume that 
$p \geq 3$.

\proclaim {Lemma 1} Suppose $\mu \leq n$ and
$$s_j = -1 \quad(\text {\rm mod\,\,} p)\,, \qquad j=1,2,\ldots,\mu-1\,.$$
Then
$$a_{up+1} = a_{up+2} = \cdots = a_{up+r}\,,$$
whenever $1 \leq up+r \leq \mu$, $r=1,2,\ldots, p$.
\endproclaim

\demo {Proof of Lemma 1}
We prove the lemma by induction on $\mu$.
The statement is obviously true for $\mu=1$.
Assume that the statement is true for $\mu$. 
We may assume that $\mu \neq up+1$, otherwise 
there is nothing to prove in the inductive step.
The inductive step follows from the Newton's identity 
$$\sum_{j=1}^{\mu}{a_js_{\mu+1-j}} + \mu a_{\mu+1} = 0\,,$$
which, together with the inductive assumption and the assumption 
$$s_j = -1 \quad(\text {\rm mod\,\,} p)\,, \qquad j=1,2,\ldots,\mu\,,$$
yields 
$$\sum_{j=1}^{\mu}{a_js_{\mu+1-j}} + \mu a_{\mu+1} = 
- \sum_{j=1}^{\mu}{a_j} + \mu a_{\mu+1} =
-\mu a_{\mu} + \mu a_{\mu+1} \quad(\text {\rm mod\,\,} p)\,.$$
We conclude that
$$a_{\mu} = a_{\mu+1} \quad(\text {\rm mod\,\,} p)\,,$$
and hence the lemma is true for $\mu+1$.
\qed \enddemo

\proclaim {Lemma 2} We have $s_{\mu} = -1 \quad(\text {\rm mod\,\,} p)$
for all $\mu=1,2,\ldots,2m-1-p$.
\endproclaim

\demo {Proof of Lemma 2}
We prove the lemma by induction on $\mu$. 
Assume that the lemma is true up to an even $\mu$ and we prove it for $\mu+1 \leq 2m-1-p$.  
Note that we already know that the lemma is true for even values of $\mu \leq 2m-1-p$.
Newton's identities give
$$\sum_{j=1}^{\mu+p}{a_js_{\mu+p+1-j}} + (\mu+p) a_{\mu+p+1} = 0 \tag 5$$
and
$$\sum_{j=1}^{\mu+p+1}{a_js_{\mu+p+2-j}} + (\mu+p+1) a_{\mu+p+2} = 0\,. \tag 6$$
Taking the difference of (6) and (5), we obtain
$$\split & \sum_{j=1}^{p}{a_j(s_{\mu+p+2-j}-s_{\mu+p+1-j})} + 
\sum_{j=p+1}^{\mu+p}{a_j(s_{\mu+p+2-j}-s_{\mu+p+1-j})} + a_{\mu+p+1}s_1 \cr 
+ & (\mu+p+1)a_{\mu+p+2} - (\mu+p)a_{\mu+p+1} = 0\,. \cr \endsplit \tag 7$$
By the inductive hypothesis we have
$$\sum_{j=p+2}^{\mu+p}{a_j(s_{\mu+p+2-j}-s_{\mu+p+1-j})} = 0
\quad(\text {\rm mod\,\,} p)\,. \tag 8$$
Also, one of the Newton's identities gives $a_1s_1 + a_2 = 0$, hence $s_1=-1$.

\noindent {\bf Case 1.} Suppose that $\mu+1 \leq 2m-1-p$ and $\mu \neq vp-1$.
Then (4) together with the inductive hypothesis and Lemma 1 yields that
$$a_{\mu+p+1}a_{\mu+p+2} = a_{(\mu+p+1)/2}a_{(\mu+p+1)/2+1}=1\,.$$
Combining this with (7),(8), $s_1=-1$, and $a_{p+1}=-1$, we obtain
$$s_{\mu+p+1}-s_{\mu+1}-s_{\mu+1} -(-1)s_{\mu} - a_{\mu+p+1} + a_{\mu+p+1} = 0
\quad(\text {\rm mod\,\,} p)\,. \tag 9$$
Now a crucial observation is that 
$$s_{\mu+p+1} = s_{\mu} = -1\,, \tag 10$$ 
as $\mu$ and $\mu+p+1 \leq 2m-2$ are even. 
Together with (9) this implies
$$-1 -2s_{\mu+1} -1 = 0 \quad(\text {\rm mod\,\,} p)\,,$$
and hence $s_{\mu+1} = -1 \enskip(\text {\rm mod\,\,} p) $. 
This is the statement of the lemma for $\mu+1$. 

\noindent {\bf Case 2.} Suppose that $\mu+1 \leq 2m-1-p$ and $\mu = vp-1$.
Then $s_1=-1$ gives
$$a_{\mu+p+1}s_1 + (\mu+p+1)a_{\mu+p+2} - (\mu+p)a_{\mu+p+1} = 0
\quad(\text {\rm mod\,\,} p)\,.$$
Combining this with (7),(8), and $a_{p+1}=-1$, we obtain
$$s_{\mu+p+1}-s_{\mu+1}-s_{\mu+1} -(-1)s_{\mu} = 0 \quad(\text {\rm mod\,\,} p)\,. \tag 11$$
Observe that (10) is valid in this case as well, hence it follows from (11) that
$$-1-2s_{\mu+1}-1 = 0 \quad(\text {\rm mod\,\,} p)\,,$$
and hence $s_{\mu+1} = -1 \enskip(\text {\rm mod\,\,} p)$.
This is the statement of the lemma for $\mu+1$.
\qed \enddemo

Lemmas 1 and 2 imply the following.

\proclaim {Lemma 3} We have 
$$a_{up+1} = a_{up+2} = \cdots = a_{up+r}$$
whenever $1 \leq up+r \leq n-p-1$, $r=1,2,\ldots, p$.
\endproclaim

Now we are ready to prove Theorem 2.

\demo{Proof of Theorem 2}
Suppose $P_{2m} \in {\Cal L}_{2m}$ is of the form
$$P_{2m}(z) = \sum_{j=1}^{2m+1}{a_jz^{2m+1-j}}\,, \qquad a_j \in \{-1,1\}\,,
\quad j=1,2,\ldots,2m+1\,,$$
and assume that
$$Q_{4m}(z) = P_{2m}(z)P_{2m}^*(z)\,.$$
As we have observed before, (2) implies that $P_{2m}^*(z)=(-1)^mP_{2m}(-z)$, 
and without loss of generality we may assume that $a_1=1$.
As in Lemmas 1, 2, and 3, let
$$a_1 = a_2 = \cdots = a_p =1, \quad a_{p+1} = -1\,.$$
Then (4) implies that $p$ is odd and, considering $P_{2m}^*$ rather than
$P_{2m}$ if necessary, without loss of generality we may assume that $p \geq 3$.
Let
$$n=2m+1=up+r\,, \qquad 0 \leq r \leq p-1\,.$$

\noindent {\bf Case 1.} Suppose $u \geq 3$ and $p \geq 5$. Then by Lemma 3 
we can deduce that
$$a_{p+2}= a_{p+3}=a_{p+4}\,.$$
On the other hand (2) implies that
$$a_{n-p-1}=-a_{n-p-2}=a_{n-p-3}\,.$$ 
However, this is impossible by Lemma 3.

\noindent {\bf Case 2.} Suppose $u \geq 3$, $p=3$, and $n \geq 13$.  Then by 
Lemma 3 we can deduce that
$$a_{p+4}= a_{p+5}=a_{p+6}\,.$$
On the other hand (2) implies that 
$$a_{n-p-3}=-a_{n-p-4}=a_{n-p-5}\,.$$ 
However, this is impossible by Lemma 3.

\noindent {\bf Case 3.} Suppose $u=2$, $p \geq 5$, and $r \geq 4$.
Then by Lemma 3 we have
$$a_{p+2}=a_{p+3}=a_{p+4}=-1\,.$$ 
On the other hand (2) implies that 
$$a_{n-p-1}=-a_{n-p-2}=a_{n-p-3}\,.$$
However, this is impossible by Lemma 3.

\noindent {\bf Case 4.} Suppose $u=2$, $p \geq 7$, and $r \leq 3$.
First we observe that Lemma 3 and (2) imply 
$P_{2m}(1) \geq p-4 \geq 3$. So $Q_{4m}(1) = P_{2m}(1)^2 \geq 9$, 
hence $m$ must be even.  
But then $P_{2m}(1) \geq p$. Hence
$$4p+5 = 2up+5 \geq 2n-1 = 4m+1 =Q_{4m}(1) = P_{2m}(1)^2 \geq p^2\,,$$
that is $p^2-4p-5 \leq 0$, which is impossible for $p \geq 7$.

\noindent {\bf Case 5.} Suppose $u=1$. Then (2) implies that $P_{2m}(1) \geq p-1$.
But then
$$4p-3 \geq 2(p+r)-1 \geq 2n-1 = 4m+1 \geq Q_{4m}(1) = P_{2m}(1)^2 \geq (p-1)^2\,,$$
that is,
$p^2- 6p+4 \leq 0$, which is impossible for $p \geq 7$.

The impossibility of the above five cases shows that if $n \geq 13$, 
then we must have either $u=2$ with $p \leq 5$ and $r \leq 3$, 
or else, $u=1$ with $p \leq 5$. Since $n = up+r$, both of the above cases 
is impossible. Hence $n \leq 12$ must be the case, and the theorem is proved.
\qed \enddemo

\enddocument